\documentclass[leqno]{article}
\usepackage{amsmath}
\usepackage{amssymb}

\begin{document}
\hyphenation{semi-per-fect Grothen-dieck}
\newtheorem{Lemma}{Lemma}[section]
\newtheorem{Th}[Lemma]{Theorem}
\newtheorem{Prop}[Lemma]{Proposition}
\newtheorem{OP}[Lemma]{Open Problem}
\newtheorem{Cor}[Lemma]{Corollary}
\newtheorem{Def}[Lemma]{Definition}
\newtheorem{Ex}[Lemma]{Example}
\newtheorem{Exs}[Lemma]{Examples}
\newenvironment{Remarks}{\noindent {\bf Remarks.}\ }{}
\newtheorem{Remark}[Lemma]{Remark}
\newenvironment{Proof}{{\sc Proof.}\ }{~\rule{1ex}{1ex}\vspace{0.5truecm}}
\newcommand{\End}{\mbox{\rm End}}
\newcommand{\K}{\mbox{\rm K.dim}}
\newcommand{\Hom}{\mbox{\rm Hom}}
\newcommand{\Ext}{\mbox{\rm Ext}}
\newcommand{\supp}{\mbox{\rm supp}}
\newcommand{\Max}{\mbox{\rm Max}}
\newcommand{\cl}{\mbox{\rm cl}}
\newcommand{\add}{\mbox{\rm add}}
\newcommand{\Inv}{\mbox{\rm Inv}}
\newcommand{\rk}{\mbox{\rm rk}}
\newcommand{\Tr}{\mbox{\rm Tr}}
\newcommand{\Sat}{\mbox{\bf Sat}}
\newcommand{\card}{\mbox{\rm card}}
\newcommand{\Ann}{\mbox{\rm Ann}}
\newcommand{\proj}{\mbox{\rm proj-}}
\newcommand{\codim}{\mbox{\rm codim}}
\newcommand{\B}{\mathcal{B}}
\newcommand{\Cong}{\mbox{\rm Cong}}
\newcommand{\Spec}{\mbox{\rm Spec-}}
\newcommand{\coker}{\mbox{\rm coker}}
\newcommand{\Cl}{\mbox{\rm Cl}}
\newcommand{\Ses}{\mbox{\rm Ses}}
\newcommand{\im}{\mbox{\rm Im}}
\newcommand{\Cal}[1]{{\cal #1}}
\newcommand{\+}{\oplus}
\newcommand{\N}{\mathbb N}
\newcommand{\Z}{\mathbb{Z}}
\newcommand{\C}{\mathbb{C}}
\newcommand{\T}{\mathbb{T}}
\newcommand{\R}{\mathbb{R}}
\newcommand{\notsim}{\sim\makebox[0\width][r]{$\slash\;\,$}}
\newcommand{\Mod}{\mbox{\rm Mod-}}
\newcommand{\lmod}{\mbox{\rm -mod}}
\newcommand{\mspec}{\mbox{\rm Max}}

\title{Local Morphisms and Modules with a Semilocal Endomorphism Ring}

\author{Alberto Facchini\thanks{Partially
supported by Gruppo Nazionale Strutture Algebriche e Geometriche e
loro Applicazioni of Istituto Nazionale di Alta Matematica, Italy,
and by Universit\`a di Padova (Progetto di Ateneo CDPA048343
``Decomposition and tilting theory in modules, derived and cluster
categories'').}\\
Dipartimento di Matematica Pura e Applicata, \\ Universit\`a di
Padova, 35131 Padova, Italy \and Dolors
Herbera\protect\thanks{Partially supported by the DGI and the
European Regional Development Fund, jointly, through Project
BFM2002-01390, and by the Comissionat per Universitats i Recerca
of the Generalitat de Ca\-ta\-lunya. \protect\newline 2000
}
\\ Departament de Matem\`atiques, \\
Universitat Aut\`onoma de Barcelona, \\ 08193 Bellaterra
(Barcelona), Spain}
\date{\phantom{Ciao}}

\maketitle

\section{Introduction}

Let $R$ and $S$ be rings. A ring morphism $\varphi \colon R \to S$
is said to be {\em local} if, for every $r\in R$, $r$ is
invertible in $R$ whenever $\varphi (r)$ is invertible in $S$
\cite{campsdicks}. For instance, if $R$ is a ring and $I$ is a
two-sided ideal of $R$ contained in the Jacobson radical of $R$,
the canonical projection $R \to R/I$ is a local morphism.
Conversely, the kernel of every local morphism $R \to S$ is
contained in the Jacobson radical of $R$
\cite[Lemma~3.1]{dolors2}. We will denote by $J(R)$ the Jacobson
radical of any ring $R$.

In
Algebraic Geometry and Commutative Algebra, local morphisms are defined as the ring morphisms
$\varphi\colon R\to S$, between local commutative rings
$(R,\mathcal{M})$ and $(S,\mathcal{N})$, for which $\varphi
(\mathcal{M})\subseteq \mathcal{N}$. This definition coincides with ours in the case of $R$ and $S$ local.

In this spirit, Cohn \cite{Cohnn} considered local morphisms $R
\to S$ when $R$, $S$ are not necessarily commutative and $S$ is a
division ring. It is easily seen that if a ring $R$ has a local
morphism into a division ring, then $R$ is a local ring.

Recall that a ring
$R$ is called \emph{semilocal} if
$R/J(R)$ is a semisimple artinian ring. The aim of this paper is to
prove that under weak
finiteness assumptions on an object $A$ of a Grothendieck category $\Cal
C$, the endomorphism ring
$\End_{\Cal C}(A)$ of $A$ is semilocal. We prove that these rings
$\End_{\Cal C}(A)$ are semilocal
making use of suitable ring homomorphisms which we show to be local morphisms.

It is known that endomorphism rings of artinian modules over an
arbitrary ring \cite[Corollary~6]{campsdicks}, or of finitely
generated modules over a semilocal commutative ring, or of
finite-rank torsion-free modules over a commutative valuation
domain or a semilocal commutative principal ideal domain
\cite[Lemma 2.3, Theorems~5.2 and 5.4]{warfield} are semilocal. A
number of other examples of modules with a semilocal endomorphism
ring are given in \cite{herberasham}. In this paper, we extend
these results from the category $\Mod R$ of right $R$-modules to
an arbitrary Grothendieck category $\Cal C$, and, also in the case
in which $\Cal C=\Mod R$, we obtain new classes of modules whose
endomorphism rings are semilocal. The advantage of knowing that a
module has a semilocal endomorphism ring lies in the fact that
modules with a semilocal endomorphism ring have a very good
behavior as far as direct sums are concerned: they cancel from
direct sums, satisfy the $n$-th root property, have only finitely
many direct summands up to isomorphism, and have only finitely
many direct-sum decompositions up to isomorphisms in the sense of
the Krull-Schmidt theorem \cite[\S~4.2]{libro}. Moreover,  classes
of modules with semilocal endomorphism rings give rise to Krull
monoids \cite[Theorem~3.4]{JA2002}. This implies that though
modules with semilocal endomorphism rings do not have uniqueness
of direct-sum decomposition up to isomorphism, the direct-sum
decompositions of these modules have a very regular geometrical
pattern. Conversely, every finitely generated Krull monoid arises
in this way from a finitely generated module over a noetherian
commutative semilocal ring \cite{Wiegand}.

After a first introductory section with the main elementary
properties of local morphisms
(Section~\ref{Local morphisms}), we prove in Section~\ref{j} that every
finitely presented module over a semilocal ring has a semilocal
endomorphism ring
(Theorem~\ref{srt}). This is one of the main results of the paper,
and generalizes the previously
known fact that every  finitely generated module over a commutative
semilocal ring has a semilocal
endomorphism ring. We give an example of a finitely generated module
over a noncommutative
semilocal ring whose endomorphism ring is not semilocal
(Example~\ref{nosemilocal}).

In Section~\ref{sspectral}, we show that local morphisms arise
naturally in the construction of the spectral category $\Spec\Cal
C$ of an arbitrary Grothendieck category $\Cal C$. The spectral
category is obtained from $\Cal C$ inverting all essential
monomorphisms \cite{gabrieloberst}, and there is a natural functor
$P\colon\Cal C\to\Spec\Cal C$. If $A$ is an object of $\Cal C$,
there is a close relation between the fact that the ring morphism
$\varphi _A\colon \End _{\mathcal{C}}(A) \to \End _{\Spec
\mathcal{C}}(A)$ induced by the functor $P$ is local and that fact
that every monomorphism $A\to A$ is an isomorphism. This allows us
to generalize \cite[Theorem~3(1)]{herberasham}. In particular, a
corollary of this is that endomorphism rings of artinian modules
are semilocal.

In Section~\ref{Finitely copresented objects}, we consider {\em
finitely copresented} objects, that is,
the objects $A$ of a Grothendieck category $\Cal C$ for which there
exists an exact sequence
$0\to A\to L_0 \to L_1\to 0$
with $L_0$ injective and both $L_0$ and $L_1$ of finite Goldie
dimension. For a finitely copresented object $A$, there is a local morphism
$\End _\mathcal{C}(A)\to \End _{\Spec\mathcal{C}}(A)\times \End
_{\Spec\mathcal{C}}(L_1)$
(Theorem~\ref{dlocal}). As a corollary, the endomorphism ring of a
finitely copresented object is
semilocal. For instance,  this shows that finite-rank torsion-free
modules over any semilocal
commutative noetherian domain $R$ of Krull dimension $1$ have
semilocal endomorphism rings
(Corollary~\ref{ext}), a fact which was previously known only under
the stronger condition of
$R$ semilocal commutative principal ideal domain.

In Section~\ref{dc}, we dualize the construction of spectral
category, obtaining a category $\Cal C'$ inverting all superfluous
epimorphisms of a Grothendieck category  $\Cal C$. There is a
natural functor $F\colon\Cal C\to\Cal C'$. As the dual of a
Grothendieck category $\mathcal{C}$ is not a spectral category in
general, the additive category $\Cal C'$ obtained in this way is
not necessarily a Grothendieck category. We consider the ring
morphism $\psi _A\colon \End _{\mathcal{C}}(A) \to \End _{\Cal
C'}(A)$ induced by the functor $F$ for every object $A$ of $\Cal
C$. This morphism $\psi_A$ is local when $A$ has finite dual
Goldie dimension and every epimorphism $A\to A$ in $\mathcal{C}$
is an isomorphism (Proposition~\ref{dualgol}). For an arbitrary
object $A$ of $\mathcal{C}$, the ring morphism $(\varphi _A,\psi
_A)\colon \End _{\mathcal{C}}(A)\to \End _{\Spec
\mathcal{C}}(A)\times \End _{\mathcal{C}'}(A)$ turns out to be
local. This leads to a generalization of \cite[Theorem~3, (2) and
(3)]{herberasham}.

In the last section, we apply the
results about the functor $F$ obtained in Section~\ref{dc}  to
objects with a projective cover. For
every exact sequence $0\to K\to P\to A\to 0$ where $P\to A$ is a
projective cover, there is a
local morphism
$\End _\mathcal{C}(A)\to \End _{\mathcal{C}'}(A)\times \End
_{\mathcal{C}'}(K)$ (Theorem~\ref{dlocal'}).
Thus if both $A$ and
$K$ have finite dual Goldie dimension, then the endomorphism ring of
$A$ is semilocal.

Our rings are associative and have an identity, and modules are
unital.

\section{Local morphisms}\label{Local morphisms}

In the next lemma we collect some  basic properties of local
morphisms. If $\varphi \colon R\to S$ is a ring morphism, we shall
denote by $M_n(\varphi)\colon
M_n(R)\to M_n(S)$ the ring morphism induced by $\varphi$ between the
rings of $n\times n$ matrices with
entries in $R$ and $S$ respectively.

\begin{Lemma}\label{basic} Let $\varphi \colon R\to S$,  $\psi \colon
S\to T$ be ring
morphisms.

\noindent\emph{(1)} If $\varphi $ is local, then $\ker (\varphi)\subseteq J(R)$
{\em \cite[Lemma~3.1]{dolors2}}.

\noindent\emph{(2) }   If $\varphi
  $ is onto and local, then $\varphi
(J(R))=J(S)$ and the induced morphism
$M_n(\varphi)\colon M_n(R)\to M_n(S)$ is local for every $n>1$ {\em
\cite[Lemma~3.1]{dolors2}}.

\noindent\emph{(3)} If $\varphi$ and $\psi$ are local morphisms, then
$\psi \circ \varphi$ is local.

\noindent\emph{(4)} If $\psi \circ \varphi$ is a local morphism, then
$\varphi$ is local.
\end{Lemma}

Local morphisms can be characterized in terms of endomorphisms
between cyclic projective modules:

\begin{Lemma} \label{cyclicproj}
Let $\varphi \colon R \to S$  be a ring morphism. The
following statements are equivalent:

\noindent\emph{(1)}
The morphism $\varphi $ is local.

\noindent\emph{(2)} If $f\colon P\to P$ is an endomorphism of a cyclic
projective right $R$-module $P$ such that $f\otimes _RS$ is
invertible, then $f$ is invertible.

\noindent\emph{(3)} If $g\colon Q\to Q$ is an endomorphism of a cyclic
projective left $R$-module $Q$ such that $S\otimes _R g$ is
invertible, then $g$ is invertible.\end{Lemma}

Most of our examples of local morphisms will satisfy stronger
properties also. To avoid confusion, it is interesting to keep in mind
the following examples.

\begin{Exs}\label{dolent}
{\rm (1) For any ring $R$, the canonical projection $\pi \colon R\to
R/J(R)$ is a local morphism.

(2) If $D$ is a division ring, the ring embedding $\varphi
\colon R=\begin{pmatrix} D& D\\0&D
\end{pmatrix} \to S=M_2(D)$ is local. Let $e_1=\begin{pmatrix} 1& 0\\0&0
\end{pmatrix}$,   $e_2=\begin{pmatrix} 0& 0\\0&1
\end{pmatrix}$ and $x=\begin{pmatrix} 0& 1\\0&0
\end{pmatrix}$. Left multiplication by $x$ induces a
morphism $f\colon e_2R\to e_1R$ such that $M_n(f\otimes _RS)$
is invertible for all $n\ge 1$, but $f$ is not invertible.
Notice that, in view of (1), another local morphism of $R$ is given by the
natural projection $R\to R/J(R)\cong D\times D$.

(3) It is not true in general that $\varphi \colon R \to S$ local
implies that the induced morphism between the matrix ring
$M_n(\varphi )\colon M_n(R)\to M_n(S)$ is local for all $n\ge 2$.
An example in which this fails is given in \cite[p.~189]{dolors2}}.
\end{Exs}

In Section \ref{dc}, we shall recall the definition of the dual
Goldie dimension $\codim(A)$ of an object $A$ of an arbitrary
Grothendieck category. It is a non-negative integer or $\infty$. Now
we only recall that a ring $R$
is semilocal if and only if $\codim(R_R)$ is finite, if and only if
$\codim(_RR)$ is finite. In this case, $\codim(R_R)=\codim(_RR)$
is the Goldie dimension $\dim(R/J(R))$ of the semisimple module $R/J(R)$
(cf.~\cite[Proposition 2.43]{libro}).

The following  deep result by Rosa Camps and Warren Dicks (see
\cite[Theorem~1]{campsdicks} or \cite[Theorem~4.2]{libro})
characterizes semilocal rings in terms of local morphisms. We will
use it throughout the paper.

\begin{Th}\label{semilocal'} If $R\to S$ is a local morphism between
arbitrary rings $R$ and $S$, then
$\codim(R)\le \codim(S)$. In particular, a ring $R$ is semilocal if
and only if there exists a local
morphism of $R$ into a semilocal
ring,  if and only if there exists a local
morphism of $R$ into a semisimple artinian
ring.\end{Th}

In general, little can be said about rings having local
morphisms to arbitrary products of division rings or to
products of rings of matrices over division rings. If $R$ is a commutative
ring with maximal spectrum $\mspec(R)$, then the morphism $R\to
\displaystyle{\prod _{\mathcal{M}\in \mspec(R)}R/\mathcal{M}}$
given by $r\mapsto (r+\mathcal{M})$ is local.  This result can be
extended to the   noncommutative setting taking as spectrum the
set of primitive ideals of $R$.
In some sense, the Camps-Dicks Theorem characterizes semilocal
rings as those having  ``finite spectrum". Notice that from Theorem~\ref{dlocal} it will follow that for every ring $R$ there exists a local morphism of $R$ into a von
Neumann regular right self-injective ring.

If $\varphi \colon R\to S$ is a local morphism with $S$ semilocal,
it is not clear which relation there is between $R/J(R)$ and
$S/J(S)$, apart from the fact that $\codim(R)\le \codim(S)$, that
is, $\dim(R/J(R))\le \dim(S/J(S))$ (cf.~Example~\ref{dolent} and
Theorem~\ref{semilocal'}). In the following Proposition, whose
proof is modelled by the proof of \cite[Lemma~3.2]{campsmenal}, we
analyze the case in which $S/J(S)$ is a finite direct product of
division rings. We  show   that the induced morphisms
$M_n(\varphi)\colon M_n(R)\to M_n(S)$ are also local, which is not
true for arbitrary rings \cite[p.~189]{dolors2}.

\begin{Prop} \label{producte} Let   $\varphi \colon R\to S$ be a
local morphism. Assume that
$S/J(S)\cong D_1\times \dots \times D_k$, where $D_i$ is a division
ring for every $i=1,\dots ,k$. For $i=1,\dots ,k$, let $\tau _i\colon
R\to D_i$ denote the composition of $\varphi$
with the projection
$S\to D_i$.
Then:

\noindent\emph{(1)} There exist $m\le k$ and $\{i_1,\dots ,i_m\}\subseteq
\{1,\dots ,k \}$ such that $R/J(R)\cong D'_1\times \dots \times
D'_m$, where $D'_j$ is a division subring of $
D_{i_j}$ for every $j=1,\dots ,m$.

\noindent\emph{(2)} The ring $R$ has exactly $m$ maximal ideals, and
these are the ideals $\ker
(\tau _{i_j})$ for
$j=1,\dots ,m$. Hence,
\[(\tau _{i_1},\dots ,\tau _{i_m})\colon R\to D_{i_1}\times \dots
\times D_{i_m}\]
is a local morphism with kernel $J(R)$.

\noindent\emph{(3)}  The induced ring morphism $M_n(\varphi)\colon
M_n(R)\to M_n(S)$ is local for every $n\ge 1$.
\end{Prop}

\begin{Proof} Let $\pi\colon S\to S/J(S)$ denote the canonical
projection. By Lemma~\ref{basic}(3 and 4), the morphism $\pi \circ
\varphi$ is local, and $M_n(\varphi)$ is local if and only if
$M_n(\pi \circ \varphi)$ is local for any $n\ge 1$. Thus, to prove
the Proposition, we may assume that $S= D_1\times \dots \times
D_k$.

By Lemma~\ref{basic}(1), $\ker
(\varphi)\subseteq J(R)$. The inclusion $\epsilon \colon
\varphi(R)\hookrightarrow S$ is a
local morphism. Hence, for any $n\ge 1$, the morphism $M_n(\epsilon)$
is local if and only if the morphism $M_n(\varphi)$ is local. Thus, to prove
the Proposition, we may assume that $R$ is a subring of $S=D_1\times
\dots \times D_k$ such that the embedding $\varphi \colon
R\hookrightarrow D_1\times \dots \times D_k$ is local.

If $k=1$, $R$ is a division subring of $D_1$ and (1), (2), (3) hold
trivially. Now we shall proceed by
induction on
$k$. Assume $k>1$.

If $R$ has a nontrivial idempotent $e$, then $e$ is central
because all idempotents of $D_1\times \dots \times D_k$ are
central. Therefore there is a partition of $\{1,\dots ,k \}$ into
two nonempty subsets $I$, $J$ such that the embeddings
\[\varphi|_{eR}\colon eR\to eS=\prod _{i\in I}D_i \mbox{ and }
\varphi|_{(1-e)R}\colon (1-e)R\to
(1-e)S=\prod _{j\in J}D_j\] are local morphisms. By the inductive
hypothesis, the Proposition
holds for $eR$ and $(1-e)R$, so it holds for $R=eR\times (1-e)R$.
Therefore we may assume that $R$ has no nontrivial idempotent.

For any element $r\in R$, set \[ \mathrm{supp}(r)=\{\,i\mid i=1,\dots
,k,\ \tau _i
(r)\neq 0\,\}.\] Let $d$ be the function of $R$ into the set of
nonnegative integers defined by
$d(r)=\vert
\mathrm{supp} (r)\vert$, so that $d$ is a nonzero function. Let
$\ell$ be the least nonzero
value of $d$. If $\ell =k$, then we can choose any
$i\in \{1,\dots ,k\}$ and we get that $\tau _i \colon R\to D_i$ is
local, so the Proposition follows from the case $k=1$. Assume
$\ell <k$.

Let $r\in R$ be an element such that $d(r)=\ell$. Suppose that there
exists $t\in R$
such that $1-tr$ is not invertible. Then $tr\ne 0$, so that
$\mathrm{supp}(tr)\subseteq\mathrm{supp}(r)$ implies
$\mathrm{supp}(tr)=\mathrm{supp}(r)$. As $\mathrm{supp}(r) \cup
\mathrm{supp}(1-tr)=\{1,2,\dots, k\}$ and $1-tr$ is not invertible, that is,
$\mathrm{supp}(1-tr)\ne\{1,2,\dots, k\}$, it follows that
$\mathrm{supp}(r)\not\subseteq
\mathrm{supp}(1-tr)
$. This implies that $\mathrm{supp}( r(1-tr))= \mathrm{supp}(r)\cap
\mathrm{supp}(1-tr)\subsetneq
\mathrm{supp}(r)$. Hence, by the choice of $r$, $r=rtr$. But then
$tr$ is idempotent and, as $R$ has no nontrivial idempotent and
$tr\ne 0$, it follows that
$tr=1$, which is impossible because $d(r)=\ell <n$. This shows that
$1-tr$ is invertible for every $t\in R$,  so that $r\in J(R)$.

Let \[ I=\{x\in R\mid \mathrm{supp}(x)\subseteq
\mathrm{supp}(r)\}= \{x\in R\mid \mathrm{supp}(x)=
\mathrm{supp}(r)\}\cup \{0\}.\]  Note that $I$ is a two-sided
ideal of $R$ which is contained in $J(R)$ by our previous
argument. Set $K=\{1,\dots ,n\}\setminus \mathrm{supp}(r)$. The
local embedding $\varphi$ induces an injective ring morphism
$\overline{\varphi}\colon R/I \to \prod _{k\in K}D_k$. Let us
prove that $\overline{\varphi}$ is local. If $t\in R$ is such that
$\overline{\varphi}(t+I)$ is invertible, then $K\subseteq
\mathrm{supp}(t)$, so that $\mathrm{supp}(t)\cup
\mathrm{supp}(r)=\{1,\dots ,n\}$. As $\mathrm{supp}(t)\cap
\mathrm{supp}(r)=\mathrm{supp}(tr)$ must have either $0$ or $\ell$
elements, it follows that either $\mathrm{supp}(t)\cap
\mathrm{supp}(r)=\emptyset$ or $\mathrm{supp}(r)\subseteq
\mathrm{supp}(t)$. If $\mathrm{supp}(t)\cap
\mathrm{supp}(r)=\emptyset$, then $\varphi(t+r)$ is invertible. If
$\mathrm{supp}(r)\subseteq \mathrm{supp}(t)$, then
$\mathrm{supp}(t)=\{1,\dots ,n\}$, so that $\varphi(t)$ is
invertible. Hence either $\varphi(t+r)$ or $\varphi(t)$ is
invertible, so that either $t+r$ or $t$ is invertible in $R$. In
both cases $t+I$ is invertible in $R/I$. Thus $\overline{\varphi}$
is a local morphism. By the inductive hypothesis, claims (1),(2)
and (3) hold for $R/I$. Therefore they hold for $R$, because
$I\subseteq J(R)$.
\end{Proof}

For further reference, we specialize Proposition~\ref{producte} to
the case $k=2$.

\begin{Cor}\label{dos} Let $\varphi \colon R\to D_1\times D_2$ be
a local morphism where $D_1$ and $D_2$ are division
rings. For $i=1,2$, let $\tau _i\colon R\to D_i$ be the
composition of $\varphi$ with the projection $D_1\times D_2 \to
D_i$. Then there are two possibilities:

\noindent\emph{(1)} either $R$ is local, and there exists $i\in
\{1,2\}$ such that
$\tau _i$ is a local morphism. In this case the maximal ideal of $R$ is
$\ker(\tau _i)$;

\noindent\emph{(2)} or $R/J(R)\cong D'_1\times D'_2$, where
$D'_i$ is a division subring of $D_i$ for $i=1,2$. Moreover, $J(R)=\ker
\varphi$ and the two maximal ideals of $R$ are $\ker(\tau _1)$ and
$\ker(\tau _2)$.\end{Cor}

We conclude this section with a result that is easy but very
useful in producing examples of modules whose endomorphism ring is
semilocal.

\begin{Prop}\label{ringext}
Let $R\to S$ be a ring morphism, and let $M_S$ be an
$S$-module with $\End(M_R)$ semilocal. Then $\End(M_S)$ is semilocal.
\end{Prop}

\begin{Proof}
Since $S$-module endomorphisms of $M_S$ are $R$-module
endomorphisms, there is an embedding $\End(M_S)\to\End(M_R)$,
which is clearly a local morphism. The Proposition follows
from Theorem~\ref{semilocal'}.
\end{Proof}

\section{Finitely presented modules over semilocal rings}\label{j}

We begin this section with a known result, of which we give an
elementary proof using the notion of local morphism studied in
this paper.

\begin{Prop}\label{vhjli} Every finitely generated module over a
commutative semilocal ring has a semilocal endomorphism
ring.\end{Prop}

\begin{Proof} Let $M_R$ be a finitely
generated module over a commutative semilocal ring $R$ and let
$\End(M_R)$ be its endomorphism ring. Consider the canonical
mapping $\varphi\colon \End(M_R)\to \End(M_R/M_RJ(R))$. This
mapping $\varphi$ is a local morphism. To see it, let $f$ be an
endomorphism of $M_R$ with $\varphi(f)$ an automorphism of
$M_R/M_RJ(R)$. By Nakayama's Lemma, $f$ must be an epimorphism.
Then $f$ must be also injective by \cite[Proposition
1.2]{Vasconcelos}. This proves that $\varphi$ is  local. But
$\End(M_R/M_RJ(R))$ is the endomorphism ring of a finitely
generated module over the ring $R/J(R)$, which is a direct product
of finitely many fields. Thus $\End(M_R/M_RJ(R))$ is semilocal, so
that $\End(M_R)$ is semilocal by
Theorem~\ref{semilocal'}.\end{Proof}

Combining Proposition
\ref{vhjli} and Proposition \ref{ringext} we get the following
extension of Proposition \ref{vhjli} \cite[Lemma 2.3]{warfield}.

\begin{Prop}\label{ycfl} Let $R$ be a semilocal commutative ring and
$S$ a (not necessarily commutative) $R$-algebra. If $M_S$ is any
$S$-module such that $M_R$ is
finitely generated, then the endomorphism ring of $M_S$ is semilocal.\end{Prop}

Now we show that Proposition \ref{vhjli} can be extended to semilocal
rings not necessarily
commutative, provided we consider finitely presented modules only.

\begin{Th}\label{srt} The endomorphism ring of a finitely presented module over
a semilocal ring is a semilocal ring.\end{Th}

\begin{Proof} Let $R$ be a semilocal ring, $M$ a finitely presented
right $R$-module and $\End(M)$ its endomorphism ring.

\smallskip

{\em Step 1. The Theorem holds under the additional hypothesis
that there exists an exact sequence $0\to K\stackrel{\displaystyle
\iota}{\longrightarrow} F\to M\to 0$, where $F$ denotes a finitely
generated free $R$-module, $K$ is a submodule of $FJ(R)$ and
$\iota\colon K\to F$ denotes the inclusion.}

For every endomorphism $g$ of a right $R$-module $A$, we shall
denote by $\overline{g}$ the endomorphism of the module $A/AJ(R)$
induced by $g$.

If $f\in \End(M)$, then there exist an endomorphism $f_0$ of $F$
and an endomorphism $f_1$ of $K$ making the diagram
$$\begin{array}{ccccccccc}
0&\to& K&\stackrel{\displaystyle
\iota}{\longrightarrow} & F&\to& M&\to& 0\\
&&\phantom{f_1}\downarrow f_1 &&\phantom{f_0}\downarrow f_0
&&\phantom{f}\downarrow f &&\\  0&\to& K&\stackrel{\displaystyle
\iota}{\longrightarrow} & F&\to& M&\to& 0
\end{array}$$ commute.

We claim that the endomorphism $\overline{f_1}$ of $K/KJ(R)$ does
not depend on the choice of the lifting $f_0$ of $f$. In order to
prove the claim, let $f'_0$ be another lifting of $f$ and $f'_1$
the corresponding restriction to $K$. Then $(f_0-f'_0)(F)\subseteq
K$, so that there exists a morphism $g\colon F\to K$ with
$f_0-f'_0=\iota g$. Thus $g(K)\subseteq g(FJ(R))\subseteq KJ(R)$.
Therefore $(f_1-f'_1)(K)=(f_0-f'_0)(K)=\iota g(K)\subseteq KJ(R)$.
That is, $\overline{f_1}-\overline{f'_1}=\overline{f_1-f'_1}=0$.
This proves the claim.

Define $\psi \colon \End(M)\to \End(M/MJ(R))\times \End (K/KJ(R))$
by $\psi (f)=(\overline{f},\overline{f_1})$ for every $f\in
\End(M)$. This is a well defined mapping by the claim, and it is
clearly a ring morphism. We shall now prove that $\psi$ is local. Let
$f\in \End(M)$ be an
endomorphism of
$M$ with $\psi(f)$ invertible. Let $f_0\in \End(F)$ be a lifting
of $f$ and $f_1\colon K\to K$ be the restriction of $f_0$ to $K$.
As $\psi(f)=(\overline{f},\overline{f_1})$ is invertible,  $f$
must be surjective by Nakayama's Lemma and $\ker f\subseteq MJ(R)$.
Similarly, as
$\overline{f_1}$ is invertible, $f_1$ must be surjective and $\ker
f_1\subseteq KJ(R)$. Since $f$ is surjective, it follows that
$f_0(F)+K=F$, hence $f_0$ also must be surjective by Nakayama's
Lemma. Thus $f_0$ must be a splitting epimorphism, because $F$ is
projective, so that $F\cong F\oplus \ker f_0$. In particular, $\ker
f_0$ is a finitely generated
$R$-module, and the finitely generated semisimple modules
$F/FJ(R)$ and
$ F/F(R)\oplus \ker f_0/\ker f_0J(R)$ are isomorphic, so that $\ker
f_0/\ker f_0J(R)=0$, from which $\ker
f_0=0$. This proves that $f_0$ is an automorphism. Since $f_1$ is
surjective, that is, $f_0(K)=K$, it
follows that
$f_0^{-1}(K)=K$. Therefore $f$ is injective. This proves that
$\psi$ is a local morphism. By Theorem~\ref{semilocal'}, the
ring $\End(M)$ is semilocal.

\smallskip

{\em Step 2. For every simple $R$-module $S$ there exists a
finitely presented $R$-module $N$ such that $S\cong N/NJ(R)$.}

As $R$ is semilocal, $S$ is isomorphic to a direct summand of
$R/J(R)$, so that there exists an isomorphism $\varphi\colon
S\oplus T\to R/J(R)$ for some $R/J(R)$-module $T$. The $R/J(R)$-module $T$
is cyclic. Let $t$ denote a generator of
$T$, and let $r\in R$ be such that $\varphi(t)=r+J(R)$. Then
$N=R/rR$ has the required property, because $$\begin{array}{l}N/NJ(R)\cong
(R/rR)/(R/rR)J(R)\cong (R/rR)/(J(R)+rR/rR)\\ \qquad \cong R/(J(R)+rR)\cong
(R/J(R))/(J(R)+rR/J(R))\\ \qquad \cong
(R/J(R))/r(R/J(R))=(R/J(R))/(r+J(R))(R/J(R))\cong S.\end{array}$$

\smallskip

{\em Step 3. For every finitely generated $R$-module $M$, there
exists a finitely presented $R$-module $N$ such that
$M/MJ(R)\oplus N/NJ(R)$ is a free $R/J(R)$-module.}

The $R$-module $M/MJ(R)$ is finitely generated and semisimple, and
$R/J(R)$ contains a direct summand isomorphic to every simple
$R$-module. Therefore $M/MJ(R)$ is isomorphic to a direct summand
of $(R/J(R))^n$ for some nonnegative integer $n$. Thus there exist
simple $R$-modules $S_1,\dots,S_m$ with $M/MJ(R)\oplus
S_1\oplus\dots\oplus S_m\cong (R/J(R))^n$. By Step 2, there exist
finitely presented $R$-modules $N_1,\dots,N_m$ with
$$M/MJ(R)\oplus N_1/N_1J(R)\oplus\dots\oplus N_m/N_mJ(R)\cong
(R/J(R))^n.$$ The module $N=N_1\oplus \dots\oplus N_m$ has the
required properties.

\smallskip

{\em Step 4. Every finitely presented $R$-module $M$ has a
semilocal endomorphism ring.}

By Step 3, there exists a finitely presented $R$-module $N$ such
that
$$M/MJ(R)\oplus N/NJ(R)\cong (R/J(R))^n$$ for some $n\ge 0$. Let $F$
be the free $R$-module $R^n$, so that there exists a surjective
morphism of $R$-modules $F\to M/MJ(R)\oplus N/NJ(R)$ with
kernel $FJ(R)$. Thus there exists a surjective morphism of
$R$-modules $F\to M\oplus N$ whose kernel $K$ is contained in
$FJ(R)$. By Step 1, the finitely presented $R$-module $M\oplus N$
has a semilocal endomorphism ring. As direct summands of modules
with semilocal endomorphism rings have semilocal endomorphism
rings \cite[Proposition 1.13]{libro}, the module $M$ also has a
semilocal endomorphism ring.\end{Proof}

\begin{Remark}{\rm We have made the proof of
Theorem~\ref{srt} as self-contained as possible, but in the rest
of the paper we will develop and refine the ideas and the
techniques we have met in the proof. Step~1 in the proof of
Theorem~\ref{srt} is a consequence of Theorem~\ref{finitelycop'},
because finitely generated modules over a semilocal ring have
finite dual Goldie dimension. The remaining part of the proof of
Theorem~\ref{srt} is devoted to showing the somewhat interesting
fact that every finitely presented module over a semilocal ring is
a direct summand of a finitely presented module with a projective
cover.}
\end{Remark}

\medskip

In Example \ref{nosemilocal} we shall show that there exist
finitely generated modules over semilocal rings whose endomorphism
ring is not semilocal. Thus Proposition \ref{vhjli} cannot be
extended to arbitrary semilocal rings, and Theorem \ref{srt}
cannot be extended to arbitrary finitely generated modules.

\medskip

   Recall that a  \emph{semiperfect
ring} is a semilocal ring whose idempotents can be lifted modulo
the Jacobson radical. A \emph{semiprimary ring}  is a semilocal
ring whose Jacobson radical is nilpotent, and a \emph{right
perfect ring} is a semilocal ring whose Jacobson radical is right
$T$-nilpotent. Bj\"{o}rk proved that finitely presented right modules
over a semiprimary ring have a semiprimary endomorphism ring
\cite[Theorem~4.1]{bjork2}. This   result was reproved and
extended by Schofield \cite[Theorem~7.18]{schofield} and
Rowen~\cite[Corollary~11]{rowen}. Their results show that a
finitely presented right module over a right (or left) perfect
ring has a right (left, respectively) perfect endomorphism ring.
Wiegand  constructed plenty of examples of finitely generated
modules over local (in particular, semiperfect) commutative
noetherian rings whose endomorphism rings are semilocal but not
semiperfect \cite{Wiegand}.

Our next example is a variation of \cite[Example~2.1,
p.~127]{bjork}. It shows that the endomorphism ring of finitely
generated modules over semiprimary rings need not be semilocal.

\begin{Ex} \label{nosemilocal}
{\rm Let $K$ be a field with a non-onto endomorphism $\alpha
\colon K\to K$. Let $K_0=\alpha (K)$. Let ${}_KV$ be a non-zero
$K$-vector space. View ${_KV}$ as a $K$-$K$-bimodule taking the
scalar product by $K$ as left action and setting as right action
$v\cdot k= \alpha(k)v$ for every $v\in V$ and every $k\in K$.

  Let $R=\begin{pmatrix} K& {}_KV_K\\0&K
\end{pmatrix}$. Then $J(R)=\begin{pmatrix} 0& {}_KV_K\\0&0
\end{pmatrix}$, $R/J(R)\cong K\times K$ and $J(R)^2=0$, so that $R$ is
semiprimary. Fix  $a\in K\setminus K_0$ and $0\neq w\in V$.
Consider the right ideal
\[I=\sum _{n\ge 0}\begin{pmatrix} 0& a^nw\\0&0
\end{pmatrix}R=\begin{pmatrix} 0& K_0[a]w\\0&0
\end{pmatrix}\] of $R$.
Then $E:=\End _R(R/I)\cong \mathcal{I}/I$, where $\mathcal{I}$ is
the idealizer of $I$ in $R$, that is, $\mathcal{I}=\{r\in R\mid
rI\subseteq I\}.$

Let $\begin{pmatrix} k_1& v\\0&k_2
\end{pmatrix}\in \mathcal{I}$. As
\[\begin{pmatrix} k_1& v\\0&k_2
\end{pmatrix}\begin{pmatrix} 0& w\\0&0
\end{pmatrix}=\begin{pmatrix} 0& k_1w\\0&0
\end{pmatrix} \in {I},\]
we deduce that $\mathcal{I}=\begin{pmatrix} K_0[a]& V\\0&K
\end{pmatrix}$. Hence $E/J(E)\cong K_0[a]\times K$.

If we choose $K$, $\alpha$ and $a$ such that $a$ is transcendental
over $K_0$, then $K_0[a]\times K$ is not semisimple artinian.
Hence, $E$ is not semilocal.}
\end{Ex}

If in Example~\ref{nosemilocal} one considers  the trivial
extension of $K$ by $V$  instead of $R$, that is,
\[K\ltimes V=\left\{\left.\begin{pmatrix} k& v\\0&k
\end{pmatrix}\right| k\in K\mbox{ and }v\in V\right\},\]
then one obtains an example of a cyclic module over the local ring $K\ltimes
V$ whose endomorphism ring is not semilocal.

\section{Spectral Categories} \label{sspectral}

In this section, we shall recall some results about spectral
categories that will be used in the following section. Spectral
categories were introduced by Gabriel and Oberst,
see \cite{gabrieloberst} or \cite[Ch.~V, \S 7]{stenstrom}. For a
Grothendieck category $\mathcal{C}$, the \emph{spectral category}
of $\mathcal{C}$, denoted by $\Spec \mathcal{C}$, is the category
with the same objects as $\mathcal{C}$ and, for objects $A$ and
$B$ of $\mathcal{C}$, with $\Hom _{\Spec
\mathcal{C}}(A,B)=\varinjlim \Hom _{\mathcal{C}}(A',B)$, where the
direct limit is taken over the downwards directed family of
essential subobjects $A'$ of $A$. There is a left exact canonical
functor $P\colon \mathcal{C}\to
\Spec \mathcal{C}$, which is the identity on objects and takes
$f\in \Hom _{\mathcal{C}}(A,B)$ to its canonical image in $\Hom
_{\Spec \mathcal{C}}(A,B)$. This functor $P$ induces a   ring
morphism
\[\varphi _A\colon \End _{\mathcal{C}}(A) \to \End _{\Spec
\mathcal{C}}(A)\] for every object $A$ of $\mathcal{C}$.

\begin{Remark}\label{spectral}{\rm
The
kernel of $\varphi _A$ is the ideal $I_A$ of all $f\in \End
_{\mathcal{C}}(A)$ with kernel essential in~$A$.}
\end{Remark}

For every object $A$ of $\mathcal{C}$, let $E(A)$ denote the
injective envelope of $A$ in $\mathcal{C}$. Then $\End
_{\Spec \mathcal{C}}(A)\cong \End
_{\mathcal{C}}(E(A))/J(\End _{\mathcal{C}}(E(A)))$ is a von
Neumann regular right self-injective ring.

\begin{Remark}\label{spectral'}{\rm
If $A$ is an injective object, the morphism $\varphi _A\colon \End _{\mathcal{C}}(A) \to \End _{\Spec
\mathcal{C}}(A)$ is the canonical projection of $\End _{\mathcal{C}}(A)$ onto $\End
_{\mathcal{C}}(A)/J(\End _{\mathcal{C}}(A))\cong \End
_{\Spec \mathcal{C}}(A)$. Therefore $\varphi _A$ is a local morphism for every injective object $A$  (Example~\ref{dolent}(1)).}
\end{Remark}

Recall that an object $A$ of a Grothendieck category
    is said to be {\em directly finite} if it is not isomorphic to a
proper direct summand of itself.

\begin{Prop}\label{local}
Let $A$ be an object in a Grothendieck category $\mathcal{C}$. If
every monomorphism  $A\to A$ is an isomorphism, then $\varphi_A\colon
\End _{\mathcal{C}}(A) \to
\End _{\Spec \mathcal{C}}(A)$ is a local morphism. Conversely, if $\varphi
_A$ is a local morphism and $E(A)$ is directly finite, then every
monomorphism  $A\to A$ is an isomorphism. \end{Prop}

\begin{Proof}
Assume that every monomorphism  $A\to A$ is an isomorphism. Let $f\in
\End _{\mathcal{C}}(A)$. If $\varphi
_A(f)$ is invertible, then any extension $\overline{f}\colon
E(A)\to E(A)$ of $f$ is a monomorphism. Thus $f$ is a monomorphism and,
hence, an isomorphism. Conversely, let $\varphi
_A$ be a local morphism and $E(A)$ directly finite. If $f\colon
A\to A$ is a monomorphism, then $f$ extends to a monomorphism
$\overline{f}\colon
E(A)\to E(A)$. As $E(A)$ is directly finite, $\overline{f}$ is an
automorphism. Thus $\varphi_A(f)$ is invertible. Since $\varphi_A$ is
local, $f$ must be an isomorphism.\end{Proof}

\begin{Prop}\label{finitegd}
The following conditions are equivalent for an object $A$ of a
Grothendieck category $\mathcal{C}$ and a nonnegative integer $n$.

\noindent\emph{(1)} $A$ has finite Goldie dimension $n$.

\noindent\emph{(2)} $P(A)$ is an object of finite length $n$ in
$\Spec \mathcal{C}$.

\noindent\emph{(3)} $\End
_{\Spec \mathcal{C}}(A)$ is a semisimple artinian ring of Goldie
dimension $n$.
\end{Prop}

\begin{Proof}
$(1)\Rightarrow (2).$ If $A$ has finite Goldie dimension $n$, then
$P(A)\cong P(E(A))$ is a semisimple object in $\Spec \mathcal{C}$ of
composition length $n$ \cite[p. 133 ]{stenstrom}.

$(2)\Rightarrow (3).$ Every object of finite length in a spectral
category is semisimple, hence it has a semisimple artinian
endomorphism ring.

$(3)\Rightarrow (1)$. Assume $\End
_{\Spec \mathcal{C}}(A)\cong \End
_{\mathcal{C}}(E(A))/J(\End _{\mathcal{C}}(E(A)))$ is semisimple
artinian. Then $\End _{\mathcal{C}}(E(A))$ is semiperfect, thus
$E(A)$ decomposes into a finite direct sum of injective
indecomposable subobjects. Therefore $E(A)$, hence $A$, has
finite Goldie dimension.~\end{Proof}

We shall denote the Goldie dimension of $A$ by $\dim(A)$. We
conclude the section with a slight generalization of
\cite[Theorem~3(1)]{herberasham}.

\begin{Cor}\label{injectiveonto} Let $A$ be an object in a
Grothendieck category
$\mathcal{C}$. Assume that $A$ has finite Goldie dimension and
that every monomorphism  $A\to A$ is an isomorphism. Then $\End
_{\mathcal{C}}(A)$ is semilocal.~\end{Cor}

\begin{Proof} By Proposition~\ref{finitegd}, the ring $\End
_{\Spec \mathcal{C}}(A)$ is semisimple artinian, and, by
Proposition~\ref{local}, $\varphi _A$ is a local morphism. The
statement follows as an application of Theorem~\ref{semilocal'}.
\end{Proof}

 From Corollary \ref{injectiveonto} and for $\Cal C=\Mod R$, $R$ any
ring, one obtains that every
artinian module has a semilocal endomorphism ring. For a different
example, let $R$ be a
commutative ring of Krull dimension $0$, that is, such that every
prime ideal is maximal. Let $M_R$
be a finitely generated module of finite Goldie dimension. Then
$\End(M_R)$ is semilocal
\cite{Vasc}.

\section{Finitely copresented objects}\label{Finitely copresented objects}

In all this section, $\mathcal{C}$ will denote a Grothendieck
category. An object $A$ of
$\mathcal{C}$ is said to be  {\em finitely copresented} if there
is an exact sequence in $\mathcal{C}$
\[0\to A\to L_0 \to L_1\to 0, \]
with $L_0$ injective, and both $L_0$ and $L_1$ of finite Goldie
dimension.

\begin{Lemma}\label{vl} The
following statements are equivalent for an object $A$ of a
Grothendieck category $\mathcal{C}$.

\noindent\emph{(1)} The object $A$ is finitely copresented.

\noindent\emph{(2)} There is
an exact sequence
\[0\to A\to E_0 \to E_1 \]
with $E_0$ and $E_1$ injective objects of finite Goldie dimension
and $A\to E_0$   an essential monomorphism.

\noindent\emph{(3)} The object $A$ is the kernel of a morphism
between injective
objects of finite Goldie dimension.\end{Lemma}

\begin{Proof} $(1)\Rightarrow (2)$ Assume there is an exact sequence
\[0\to A\to L_0 \to L_1\to 0, \]
with $L_0$ injective, $\dim(L_0)<\infty$ and $\dim(L_1)<\infty$.
Then $L_0$ has an injective envelope $E_0$ of $A$ as a direct
summand, and the sequence
\[0\to A\to E_0 \to L_1 \]
is exact. Now substitute $L_1$ by its injective envelope $E_1$.

The implications  $(2)\Rightarrow (3)\Rightarrow
(1)$ are trivial.~\end{Proof}

Following the notation introduced in the previous section, let
$P\colon \mathcal{C}\to \Spec\mathcal{C}$ denote the canonical
functor of $\mathcal{C}$ into its spectral category.

\begin{Lemma}\label{welldefined} Let $A$ be an object of $\mathcal{C}$ and let
$L_0$ be its injective envelope. Consider the exact sequence
$0\to A\to L_0 \to L_1\to 0$, so that every $f\in \End _\mathcal{C}(A)$
extends to an endomorphism $f_0$ of $L_0$, which induces an
endomorphism $f_1$ of $L_1$. Then $P(f_1)$ depends only on $f$ and not on the
choice of the extension $f_0$ of $f$.~\end{Lemma}

\begin{Proof}
Let $f'_0$ be another extension of $f$ and $f'_1$ the
corresponding endomorphism of $L_1$. Then $f_0-f'_0$ induces a
morphism $g\colon L_0/A\to L_0$. The inverse image $g^{-1}(A)$
is essential in $L_0/A$ because $A$ is essential in $L_0$.
Therefore the endomorphism $f_1-f'_1$ of $L_1$ induced by
$f_0-f'_0$ has an essential kernel. That is,
$P(f_1-f'_1)=P(f_1)-P(f'_1)=0$.~\end{Proof}

By Lemma \ref{welldefined}, for every object $A$ of $\mathcal{C}$
there is a ring morphism
\[\chi \colon \End _\mathcal{C}(A)\to \End _{\Spec\mathcal{C}}(A)\times \End
_{\Spec\mathcal{C}}(L_1)\] defined by $\chi (f)=(P(f),P(f_1))$.

\begin{Th}\label{dlocal} The ring morphism
$\chi$ is local for every object $A$ of $\mathcal{C}$.\end{Th}

\begin{Proof}
     Let $f\in \End _\mathcal{C}(A)$ and assume $\chi(f)$ invertible. Let
$f_0\in \End _\mathcal{C}(L_0)$ be an extension of $f$,
     and let $f_1\colon L_1\to L_1$ be the induced endomorphism of
$L_1$, so that we have a
commutative diagram
$$\begin{array}{ccccccccc}
0&\to& A&\to& L_0&\to& L_1&\to& 0\phantom{.}\\
&&\phantom{f}\downarrow f &&\phantom{f_0}\downarrow f_0
&&\phantom{f_1}\downarrow f_1 &&\\  0&\to& A&\to& L_0&\to&
L_1&\to& 0.
\end{array}$$ As $P(f)$ and $P(f_1)$ are invertible,
the morphisms $f$ and $f_1$ must be essential monomorphisms. But $P(A)$ is canonically isomorphic to $P(L_0)$ and $P(f)$ is an isomorphism in $\Spec\mathcal{C}$, so that $P(f_0)$ is an isomorphism in $\Spec\mathcal{C}$. By Remark~\ref{spectral'}, 
the morphism $f_0$ of $\mathcal{C}$ is an isomorphism. 
The Snake Lemma gives an exact sequence
\[0=\ker f_1\to \mathrm{coker}f \to  \mathrm{coker} f_0=0,\]
so that $\mathrm{coker}f =0$, i.e., $f$ is also
an epimorphism.~\end{Proof}

>From Theorem~\ref{dlocal}, it follows that for every ring $R$ there exists a local morphism of $R$ into a von
Neumann regular right self-injective ring.

\begin{Th} \label{finitelycop}
Let $A$ be a finitely copresented object of a Grothendieck
category $\mathcal{C}$. Then $\End _\mathcal{C}(A)$ is a semilocal
ring. Moreover, if $L_0$ denotes the injective envelope of $A$,
then $\codim(\End
_\mathcal{C}(A))\le\dim(A)+\dim(L_0/A)$.~\end{Th}

\begin{Proof} By Proposition~\ref{finitegd}, the ring $\End
_{\Spec\mathcal{C}}(A)\times \End _{\Spec\mathcal{C}}(L_0/A)$ is
semisimple artinian and its Goldie dimension is
$\dim(A)+\dim(L_0/A)$. Now apply Theorem~\ref{dlocal} and
Theorem~\ref{semilocal'} to conclude.~\end{Proof}

In the case in which $\Cal C=\Mod R$, Theorem~\ref{finitelycop} becomes

\begin{Cor} \label{kernels}
Let $M$ be a finitely copresented right module over an arbitrary
ring $R$. Then $\End _R(M)$ is a semilocal ring. Moreover, if
$L_0$ denotes the injective envelope of $M$, then $\codim(\End
_R(M))\le\dim(M)+\dim(L_0/M)$.~\end{Cor}

We say that a module $M$ is {\em quotient finite dimensional} if
every homomorphic image of $M$ has
finite Goldie dimension.

\begin{Lemma} Let $N$ be a submodule of a module $M$. If both $N$ and
$M/N$ are quotient finite
dimensional, then $M$ is quotient finite dimensional.\end{Lemma}

\begin{Proof} Let $P$ be a submodule of $M$. We must show that $M/P$ has
finite Goldie dimension. As $M/(N+P)$ has
finite Goldie dimension, there exist injective modules $E_1,\dots,
E_n$ of Goldie dimension $1$ and a
homomorphism $f\colon M\to E_1\oplus\dots\oplus E_n$ with $\ker
f=N+P$. Now $N+P/P\cong N/N\cap P$ has
finite Goldie dimension. Hence there exist injective modules
$E_{n+1},\dots, E_m$ of Goldie dimension $1$
and a homomorphism $g\colon N+P\to E_{n+1}\oplus\dots\oplus E_m$ with
$\ker g=P$. The homomorphism $g$
extends to a homomorphism $h\colon M\to E_{n+1}\oplus\dots\oplus
E_m$. Notice that $(N+P)\cap \ker h=P$.
Consider the homomorphism $(f,h)\colon M\to E_1\oplus\dots\oplus
E_n\oplus E_{n+1}\oplus\dots\oplus E_m$.
Then $\ker(f,h)=\ker f\cap \ker h=(N+P)\cap \ker h=P$. Therefore $M/P$ has
finite Goldie dimension. \end{Proof}

\begin{Cor}\label{vjhl} A direct sum of finitely many quotient finite
dimensional modules is quotient finite dimensional.\end{Cor}

 From Lemma \ref{vl} and Corollary \ref{kernels}, we obtain:

\begin{Cor}\label{vjhl'} Every submodule of a quotient finite
dimensional injective module has a semilocal
endomorphism ring.\end{Cor}

For instance, let $R$ be a commutative noetherian semilocal domain
of Krull dimension $1$, and let $Q$ be the field of fractions of
$R$. By \cite[Theorem 1 p.~571]{matlis}, the $R$-module $Q/R$ is
artinian, so that $Q$ is a quotient finite dimensional injective
$R$-module. By Corollary~\ref{vjhl}, all $Q^n$ are quotient finite
dimensional injective modules, so that their submodules, that is,
torsion-free modules of finite rank, have semilocal endomorphism
rings. Applying Proposition~\ref{ringext}, we get the following
corollary, which generalizes a result proved by Warfield only for
the case in which $R$ is a commutative semilocal principal ideal
domain (cf.~\cite[Theorem~5.2]{warfield}).

\begin{Cor}\label{ext}
Let $R$ be a commutative noetherian semilocal domain of Krull
dimension $1$ and let $S$ be an $R$-algebra. Let $M_S$ be an
$S$-module that is torsion-free
of finite rank as an $R$-module. Then $\End(M_S)$ is semilocal.~\end{Cor}

More generally, we have shown that if $R$ is a commutative integral
domain, the field of fractions $Q$ of $R$ is a quotient finite
dimensional $R$-module,
$S$ is an $R$-algebra and $M_S$ is an $S$-module that is torsion-free
of finite rank as an $R$-module, then $M_S$ has a semilocal
endomorphism ring. For the case of $R$ a valuation domain, this is
\cite[Theorem~5.4]{warfield}.

We shall now give a further extension of \cite[Theorem~5.4]{warfield}
to the noncommutative
setting.
Recall that a right module $M$ is {\em uniserial} if its
lattice of submodules is linearly ordered by set inclusion, that is,
if for any submodules $N$ and $P$ of
$M$ either $N\subseteq P$ or $P\subseteq N$. A module is {\em serial}
if it is a direct sum of
uniserial submodules.

\begin{Cor}\label{ext'} Let $E$ be an injective serial right module
of finite Goldie dimension over an
arbitrary (not necessarily commutative) ring. Then the endomorphism
ring of every submodule of $E$ is
semilocal.~\end{Cor}

\begin{Proof} The module $E$ is a direct sum of uniserial submodules,
necessarily
finitely many because $E$ has finite Goldie dimension. Thus $E$ is
quotient finite dimensional by
Corollary~\ref{vjhl}. Now apply Corollary~\ref{vjhl'}. \end{Proof}

We conclude
this section with an application of Theorem~\ref{finitelycop} to a
category $\Cal C$ that is not a
category $\Mod R$.

\begin{Cor} \label{kernelspi} Let $R$ be a ring.  Let $E_0,E_1$ be
direct sums of $n,m$ indecomposable
pure-injective right $R$-modules, respectively. Let $f\colon E_0\to
E_1$ be a morphism whose kernel $M$ is
pure in $E_0$ and whose image $f(E_0)$ is pure in $E_1$, so that the
pure-injective envelopes of
$M$ and $f(E_0)$ are direct sums of $r\le n$ and $s\le m$ indecomposable
pure-injective right $R$-modules, respectively.  Then
$\End _R(M)$ is a semilocal ring and
$\codim(\End _R(M))\le 2r+s-n$.~\end{Cor}

\begin{Proof} Let $R\lmod$ denote the category of finitely presented
left $R$-modules, and
let $\Cal F:=\mathrm{Add}(R\lmod,\ \mbox{\rm Ab})$ denote the
category of additive functors from $R\lmod$ to the
category Ab of abelian groups. The assignment $X\mapsto X\otimes _R-$
defines a functor $\Phi \colon \Mod R\to
\Cal F$, which is a full and faithful. Moreover, $\Phi$ sends
pure-injective objects of
$\Mod R$ to injective objects of the Grothendieck category
$\Cal F$ and
pure-exact sequences of
$\Mod R$ to exact sequences of $\Cal F$
(cf.~\cite[Theorem~B.16]{jensenlenzing} or
\cite[\S~1.6]{libro}).

Therefore $\Phi$ sends the pure-exact sequences $0\to M\to E_0\to
f(E_0)\to 0$ and $0\to f(E_0)\to
E_1$ to the exact sequences $0\to \Phi(M)\to \Phi(E_0)\to
\Phi(f(E_0))\to 0$ and
$0\to \Phi(f(E_0))\to \Phi(E_1)$. Thus the sequence $0\to \Phi(M)\to
\Phi(E_0)\to
\Phi(E_1)$ is exact, i.e., the functor $\Phi(M)$ is the kernel of the
morphism $\Phi(f)\colon \Phi(E_0)\to
\Phi(E_1)$ between injective objects of $\Cal F$. Notice that if $E$
is a direct sums of $t$
indecomposable pure-injective right $R$-modules, the object $\Phi(E)$
has Goldie dimension $t$ in
$\Cal F$, because, by Proposition~\ref{finitegd}, the Goldie dimension of
$\Phi(E)$ in
$\Cal F$ is equal to the Goldie dimension of $\End_{\Spec \Cal
F}(\Phi(E))\cong \End_{\Cal
F}(\Phi(E))/J(\End_{\Cal F}(\Phi(E))\cong \End_R(E)/J(\End_R(E))$.
This shows that $\Phi(M)$ is a finitely
copresented object in
$\Cal F$ (Lemma~\ref{vl}).

Theorem~\ref{finitelycop} implies that $\End _{\Cal
F}(\Phi(M))\cong \End _R(M)$ is semilocal of dual Goldie dimension
$\le \dim(\Phi(M))+\dim(F/\Phi(M))$, where $F$ denotes the
injective envelope of $\Phi(M)$. Thus $F=\Phi(P)$, where $P$
denotes the pure-injective envelope of $M$. But
$\dim(\Phi(M))=\dim(F)=\dim(\Phi(P))=r$ and
$\dim(F/\Phi(M))=\dim(\Phi(P)/\Phi(M))=\dim(\Phi(P/M))=\dim(\Phi(E_0/M))-\dim(\Phi(E_0/P))=
\dim(\Phi(f(E_0)))-(n-r)= s-n+r$.
\end{Proof}

\section{The dual construction}\label{dc}

The construction of the spectral category can be dualized. For a
Grothendieck category $\mathcal{C}$,   consider the category
$\mathcal{C}'$ with the same objects as $\mathcal{C}$ and, for
objects $A$ and $B$ of $\mathcal{C}$, with $\Hom
_{\mathcal{C}'}(A,B)=\varinjlim \Hom _{\mathcal{C}}(A,B/B')$,
where the direct limit is taken over the upwards directed family
of superfluous (\,=\,small\,=\,inessential) subobjects $B'$ of
$B$. There is a canonical functor $F\colon \mathcal{C}\to
\mathcal{C}'$ which is the identity on objects.

More formally, assume that $\mathcal{C}$ is any abelian category
and let $S$ be the system of all its superfluous epimorphisms
(epimorphisms with superfluous
kernel), that is, the epimorphisms $s\colon A\to B$ such that, for
every subobject $A'$ of $A$, $s(A')=B$ implies $A'=A$. It is
easily seen that if $s\colon A\to B$ and $t\colon B\to C$ are
epimorphisms, then $ts$ has superfluous kernel if and only if both
$s$ and $t$ have superfluous kernels. Moreover, every co-angle
$$\begin{array}{ccc} A&\stackrel{\displaystyle g}{\longrightarrow} &A'\\
s\downarrow\phantom{s}&&\\
B&&\end{array}$$ has a pushout $$\begin{array}{ccc}
A&\stackrel{\displaystyle g}{\longrightarrow} &A'\\
s\downarrow\phantom{s}&&\phantom{s'}\downarrow s'\\
B&\stackrel{\displaystyle f}{\longrightarrow} &B',\end{array}$$
and if $s$ is a superfluous epimorphism, then $s'$ is
a superfluous epimorphism, because if $k\colon K\to
A$ is the kernel of $s$, then the kernel of $s'$ is the image of
$gk\colon K\to A'$. Thus $S$ is a left-calculable multiplicative
system of morphisms in $\mathcal{C}$ \cite[p.~152]{Popescu}.
Assume that the abelian category $\mathcal{C}$ has a set of
generators, so that it is locally small and colocally small.
Consider, for every object $B$ of $\mathcal{C}$, the category $B/S$
whose objects are the pairs $(s,C)$ with $s\colon B\to C$ a superfluous
epimorphism in $\mathcal{C}$ and whose
morphisms $f\colon (s,C)\to (s',C')$ are the morphisms $f\colon
C\to C'$ in $\mathcal{C}$ with $sf=s'$. Then $B/S$ has a small
cofinal subcategory, because it is sufficient to consider the
pairs $(s,C)$ where $s\colon B\to C$ ranges in a set of
representatives of quotient objects of $B$ in the colocally small
category $\mathcal{C}$. Under these conditions, the category
$\mathcal{C}'=\mathcal{C}_S$ of additive fractions of
$\mathcal{C}$ relative to $S$ exists
\cite[Theorem~4.1.4]{Popescu}. It has the same objects as
$\mathcal{C}$ and, for objects $A$ and $B$ of $\mathcal{C}$, $\Hom
_{\mathcal{C}'}(A,B)$ is the inductive limit of the abelian groups
$\Hom _{\mathcal{C}}(A,C)$ where $(s,C)$ ranges in $B/S$, that is,
the inductive limit of the functor $\Hom _{\mathcal{C}}(A,-)\colon
B/S\to\mbox{\rm Ab}$. The morphisms in $\mathcal{C}'$ are usually
denoted as fractions $(s/f)\colon A\to B$, where $(s,C)$ is an
object of $B/S$ and $f\colon A\to C$ is a morphism in
$\mathcal{C}$. This category $\mathcal{C}'$ can also be
constructed by passing to the dual category of $\mathcal{C}$. Let
$\mathcal{C}$ be an abelian category with a set of generators.
Then the dual category $\mathcal{C}^0$ of $\mathcal{C}$ is a
locally small abelian category, the superfluous epimorphisms of
$\mathcal{C}$ become the essential
monomorphisms in
$\mathcal{C}^0$, so that $S$ is a right-calculable system in
$\mathcal{C}^0$ \cite[Corollary~4.2.2]{Popescu} and it is possible
to construct $\mathcal{C}'=\mathcal{C}^0_S$. Notice that the
category $\mathcal{C}^0$ is locally small but does not satisfies
the hypothesis of \cite[Theorem~4.2.5]{Popescu}, so that
$\mathcal{C}'=\mathcal{C}^0_S$ is not necessarily a spectral
category.

The category $\mathcal{C}'$ defined in this way can be far from
being spectral also in the case of a Grothendieck category
$\mathcal{C}$. For instance, if $\mathcal{C}$ is the category
Ab of abelian groups, and $\Z$ is the abelian group of
integers, then $\Z$ does not have non-zero superfluous subobjects
in Ab, so that the endomorphism ring of the object
$\Z$ in the category $\mathcal{C}'$ is the ring $\Z$, while in
spectral categories endomorphism rings are always von Neumann
regular right self-injective rings. Nevertheless we are only
interested in the ring morphisms $\psi_A\colon
\End_{\mathcal{C}}(A) \to \End _{\mathcal{C}'}(A)$ induced by the
functor $F$ for every object $A$ of $\mathcal{C}$. The kernel of
$\psi_A$ is the ideal $K_A$ of all $f\in \End _{\mathcal{C}}(A)$
whose image is a superfluous subobject of $A$.

For instance, let $R$ be a ring, $\mathcal{C}=\Mod R$, $P$ a
finitely generated projective right $R$-module and $\End _R(P)$
its endomorphism ring. Then $\End _{\mathcal{C}'}(P)\cong \End
_R(P)/J(\End _R(P))$ \cite[Proposition~17.11]{AF2}. More generally,
if $N$ is a finitely generated right $R$-module with a projective
cover $P$, then $\End _{\mathcal{C}'}(N)\cong \End _R(P)/J(\End _R
(P))$. Hence, if $N/NJ(R)$ is projective as an
$R/J(R)$-module, $\End _{\mathcal{C}'}(N)\cong \End _R(N/NJ(R))$
\cite[Corollary~17.12]{AF2}.

We state the following elementary lemma for later reference.

\begin{Lemma}\label{sup} Let $f\colon A\to B$ be a morphism in a
Grothendieck category $\mathcal{C}$. Then:

\noindent\emph{(1)} The morphism $F(f)$ is an isomorphism
if and only if $f$ is a superfluous epimorphism.

\noindent\emph{(2)} If $B$ is projective and $F(f)$ is an isomorphism, then $f$ is an isomorphism.\end{Lemma}

Recall that the {\em dual Goldie dimension} $\codim(A)$ of an
object $A$ of a Grothendieck category $\mathcal{C}$ is the Goldie
dimension of the dual lattice of the lattice $\mathcal{L}(A)$ of
all subobjects of $A$ \cite[\S~2.8]{libro}. An object $A$ of
$\mathcal{C}$ is {\em couniform} if $\codim(A)=1$, that is, if it
is uniform in the dual category of $\mathcal{C}$
\cite[p.~184]{DF}. Equivalently, a non-zero object $A$ of
$\mathcal{C}$ is couniform if and only if the sum of any two
proper subobjects of $A$ is a proper subobject of $A$, if and only
if every proper subobject of $A$ is superfluous, if and only if
whenever $f\colon A'\to A$ and $g\colon A''\to A$ are morphisms in
$\mathcal{C}$ and the coproduct morphism $f\oplus g\colon A'\oplus
A''\to A$ is an epimorphism, at least one of the morphisms $f$ and $g$ is
an epimorphism. We have the following

\begin{Lemma}\label{couniform}
Let $U$ and $V$ be couniform objects of a Grothendieck category
$\mathcal{C}$. Then:

\noindent\emph{(1)}  $\End _{\mathcal{C}'}(F(U))$ is a
division ring.

\noindent\emph{(2)}  $\Hom_{\mathcal{C}'}(F(U),F(V))\ne 0$ if and only if
there exist proper subobjects $U'$ of $U$ and $V'$ of $V$ with $U/U'$
isomorphic to $V/V'$, if and only if $F(U)$ is
isomorphic to
$F(V)$.
\end{Lemma}

\begin{Proof} Every morphism $F(U)\to F(V)$ is represented by a
morphism $f\colon U\to V/V'$ for some proper subobject $V'$ of $V$.
Also, the image of such an $f\colon U\to V/V'$ is
zero in $\Hom_{\mathcal{C}'}(F(U),F(V))$ if and only if $f$ is not
an epimorphism in $\mathcal{C}$. For $U=V$, it follows that every non-zero
element of  $\End _{\mathcal{C}'}(F(U))$ is an
isomorphism, which proves (1). If $\Hom_{\mathcal{C}'}(F(U),F(V))\ne
0$, then there is an epimorphism $f\colon U\to V/V'$ for some proper
subobject $V'$ of $V$. Thus the kernel $\ker(f)\to U$ of $f$ is a
proper subobject of $U$ and $U/\ker(f)$ is isomorphic to $V/V'$. If
$U',V'$ are proper subobjects of $U,V$ respectively with $U/U'$
isomorphic to $V/V'$, then there is an epimorphism
$f\colon U\to V/V'$ and its image in $\Hom_{\mathcal{C}'}(F(U),F(V))$
is an isomorphism. The rest is clear.
\end{Proof}

\begin{Prop}\label{dualgol}
Let $A$ be an object of finite dual Goldie dimension in a
Grothendieck category $\mathcal{C}$.
Then:

\noindent\emph{(1)} The ring  $\End _{\mathcal{C}'}(F(A))$ is
semisimple artinian of
Goldie dimension $=\codim(A)$.

\noindent\emph{(2)} If $f\in \End_\mathcal{C}(A)$, the morphism
$F(f)$ is invertible if
and only if $f$ is an epimorphism.

\noindent\emph{(3)}  If every epimorphism $A\to A$ in $\mathcal{C}$ is an
isomorphism, then $\psi _A\colon \End_{\mathcal{C}}(A) \to \End
_{\mathcal{C}'}(A)$ is a local morphism.
\end{Prop}

\begin{Proof} (1) Since $\codim(A)=n$ is finite, $A$ has a
superfluous subobject $K$ with  $A/K=U_1\oplus \dots \oplus U_n$,
where $U_i$ is a couniform object for every $i=1,\dots ,n$. As $F(A)$
is isomorphic to $F(A/K)$ in $\mathcal{C}'$, we may assume that $A$
is a finite direct
sum of couniform objects $U_1,\dots , U_n$. Statement (1) is now a
consequence of Lemma \ref{couniform}.

(2) Let $f\in \End _\mathcal{C}(A)$. By Lemma~\ref{sup}(1), $F(f)$ is
invertible if and only if $f$ is a superfluous epimorphism. Since $A$
has finite dual Goldie
dimension, all epimorphisms $A\to A$ have superfluous kernels.

(3) is a consequence of (2).~\end{Proof}

Recall that in Remark~\ref{spectral} we denoted by $I_A$ the
kernel of $\varphi _A\colon \End _{\mathcal{C}}(A) \to \End
_{\Spec \mathcal{C}}(A)$, that is, the ideal of all endomorphisms
of $A$ with essential kernel, and that we denote by $K_A$ the
kernel of $\psi _A$, that is, the ideal of all endomorphisms of
$A$ with superfluous image. In   Proposition~\ref{bj} we put
together the ring morphisms $\varphi_A$ and $\psi _A$ to   obtain
a local morphism:

\begin{Prop}\label{bj}
Let $A$ be an object in a Grothendieck category $\mathcal{C}$.
Then the ring morphism \[(\varphi _A,\psi _A)\colon \End
_{\mathcal{C}}(A)\to \End _{\Spec \mathcal{C}}(A)\times \End
_{\mathcal{C}'}(A),\] defined by $f\mapsto (P(f),F(f))$ for every
$f\in \End _{\mathcal{C}}(A)$, is local. The kernel of this ring
morphism is $I_A\cap K_A$, that is, the set of all $f\in \End
_{\mathcal{C}}(A)$ with essential kernel and superfluous image.
\end{Prop}

\begin{Proof}
Assume that $f\in \End _{\mathcal{C}}(A)$ is such that $P(f)$ and
$F(f)$ are invertible.

In general, $P(f)$ is invertible if and only if $f$ is an
essential monomorphism. By Lemma~\ref{sup}(1), if $F(f)$ is
invertible, then $f$ is an epimorphism. Hence, $f$ is invertible.
\end{Proof}

\begin{Cor}
Let $A$ be an object of a Grothendieck category $\mathcal{C}$. Assume
that $A$ has
finite Goldie dimension $n$ and finite dual Goldie dimension $m$.
Then $\End _{\mathcal{C}}(A)$ is semilocal and $\codim(\End
_{\mathcal{C}}(A))\le n+m$.
\end{Cor}

\begin{Proof}
By Propositions~\ref{finitegd} and \ref{dualgol}(1), the rings $\End
_{\Spec \mathcal{C}}(A)$ and $\End
_{\mathcal{C}'}(A)$ are
semisimple artinian of Goldie dimension $\dim(A)$ and $\codim(A)$,
respectively.
Thus $\End _{\mathcal{C}}(A)$ is semilocal and $\codim(\End
_{\mathcal{C}}(A))\le n+m$ by
Theorem~\ref{semilocal'} and Proposition \ref{bj}.
\end{Proof}

If $A$ is a uniform object, then
$I_A=\{f\in \End _{\mathcal{C}}(A)\mid f\mbox{ is not a
monomorphism}\}$ and $\End _{\Spec
\mathcal{C}}(A)$ is a division ring. If $A$ is a couniform object,
then
$K_A=\{f\in \End _{\mathcal{C}}(A)\mid f\mbox{ is not an
epimorphism}\}$ and $\End _{\mathcal{C}'}(A)$ is
a division ring. Therefore, if $A$ is both uniform and couniform, the
local ring morphism $(\varphi
_A,\psi _A)$ maps $\End
_{\mathcal{C}}(A)$ into the direct product of two division rings.
 From Corollary~\ref{dos} and
Proposition~\ref{bj}, we recover the basic results on the
endomorphism ring of biuniform modules \cite[Theorem~9.1]{libro}
that we extend to the context of Grothendieck categories.

\begin{Cor} Let $A$ be an object of a Grothendieck category
$\mathcal{C}$. Assume that $A$ is uniform and couniform. Then
there are two possibilities:

\noindent\emph{(1)} either the ideals  $I_A$ and $K_A$ are comparable
and, in this case, $\End
_{\mathcal{C}}(A)$ is local with maximal ideal $I_A+K_A$, or

\noindent\emph{(2)} the ideals $I_A$ and $K_A$ are not comparable,
$\End _{\mathcal{C}}(A)/J(\End _{\mathcal{C}}(A))$ is the product
of two division rings, and $I_A$, $K_A$ are the two maximal ideals
of $\End _{\mathcal{C}}(A)$.
\end{Cor}

\section{Objects with a projective cover}

In all this section, $\mathcal{C}$ will denote a Grothendieck
category. Now we shall apply the
results of the previous section about the functor $F\colon
\mathcal{C}\to\mathcal{C}'$ to objects with a projective cover.

\begin{Lemma}\label{welldefined'} Let $A$ be an object of $\mathcal{C}$
with a projective cover $\pi\colon P\to A$, and let $K\to P$ be the
kernel of $\pi$.
Let $f\in \End _R(A)$, so that $f$ lifts to an endomorphism $f_0$ of $P$,
which restricts to an endomorphism $f_1$ of $K$. Then $F(f_1)$
depends only on $f$ and not on the
choice of the lifting $f_0$ of
$f$.~\end{Lemma}

\begin{Proof}
Let $f'_0$ be another lifting of $f$ and $f'_1$ the corresponding
restriction to $K$. As $\pi \circ (f_0-f'_0)=0$, the difference
$f_0-f'_0\colon P\to P$ factors through
the kernel $\iota\colon K\to P$ of $\pi$, that is, $f_0-f'_0=\iota g$
for a suitable morphism $g\colon
P\to K$. As
$K$ is superfluous in
$P$, its image $g(K)$ is superfluous in $K$. Therefore the image
of the restriction $f_1-f'_1\colon K\to K$ of $g$ to $K$ is
superfluous in $K$. That is,
$F(f_1-f'_1)=F(f_1)-F(f'_1)=0$.~\end{Proof}

By Lemma \ref{welldefined'}, there is a ring morphism
\[\Phi \colon \End _\mathcal{C}(A)\to \End _{\mathcal{C}'}(A)\times \End
_{\mathcal{C}'}(K)\] defined by $\Phi (f)=(F(f),F(f_1))$ for every
$f\in \End _\mathcal{C}(A)$.

\begin{Th}\label{dlocal'} Let $A$ be an object of a Grothendieck
category $\mathcal{C}$.
Suppose that there exists a projective cover $\pi\colon P\to A$. Then the ring morphism $\Phi$ is local.
\end{Th}

\begin{Proof} Let $K\to P$ be the kernel of $\pi$.
     Let $f\in \End _\mathcal{C}(A)$ be such that $\Phi(f)$ is invertible. Let
$f_0\in \End _\mathcal{C}(P)$ be a lifting of $f$,
     and let $f_1\colon K\to K$ be the restriction of $f_0$ to $K$,
     so that we have a commutative diagram
$$\begin{array}{ccccccccc}
0&\to& K&\to& P&\to& A&\to& 0\\
&&\phantom{f_1}\downarrow f_1 &&\phantom{f_0}\downarrow f_0
&&\phantom{f }\downarrow f  &&\\  0&\to& K&\to& P&\to& A&\to& 0.
\end{array}$$

     As
$F(f)$ and $F(f_1)$ are invertible, the morphisms $f$ and $f_1$ must
be epimorphisms by
Lemma~\ref{sup}(1). We must prove that $f$ is a monomorphism. As $F(P)$ and $F(A)$ are canonically isomorphic via $F(\pi)$ and $F(f)$ is an isomorphism, it follows that $F(f_0)$ is an isomorphism in $\mathcal{C}'$. From Lemma~\ref{sup}(2), we get 
that $f_0$ must be an isomorphism in $\mathcal{C}$. The Snake Lemma gives an exact sequence
$0=\ker f_0 \to \ker f \to \mathrm{coker}f_1=0.$
Hence $\ker f=0$, as we wanted to prove.\end{Proof}

\begin{Th} \label{finitelycop'}
Let $A$ be an object of a Grothendieck category $\mathcal{C}$.
Suppose that there exists a
projective cover
$\pi\colon P\to A$. Let $K\to P$ be the kernel of $\pi$ and assume
that both $A$ and
$K$ have finite dual Goldie dimension. Then $\End _\mathcal{C}(A)$
is a semilocal ring. Moreover, $\codim(\End
_\mathcal{C}(A))\le\codim(A)+\codim(K)$. \end{Th}

\begin{Proof} By
Proposition~\ref{dualgol}, the ring $\End _{\mathcal{C}'}(A)\times \End
_{\mathcal{C}'}(K)$ is semisimple artinian of dual Goldie dimension $\codim(A)+\codim(K)$. Now apply Theorems
\ref{dlocal'}  and \ref{semilocal'} to conclude.
\end{Proof}


\begin{thebibliography}{99}


\bibitem[AF]{AF2} F. W. Anderson and
K. R. Fuller, ``Rings and
Categories of Modules'', Second
Edition, Springer-Verlag, New
York, 1992.

\bibitem[B1]{bjork} J.-E.~Bj\"{o}rk,
\emph{Conditions which imply that
subrings of artinian rings are
artinian}, J. Reine Angew. Math.
\textbf{247} (1971), 123--138.

\bibitem[B2]{bjork2} J.-E.~Bj\"{o}rk,
\emph{Conditions which imply that
subrings of semiprimary rings are
semiprimary,}  J. Algebra
\textbf{19}  (1971), 384--395.

\bibitem[CD]{campsdicks} R. Camps and
W. Dicks, \emph{On semilocal
rings}, Israel J.
Math.~\textbf{81} (1993),
203--211.

\bibitem[CM]{campsmenal} R. Camps and
P. Menal, \emph{Power
cancellation for artinian
modules}, Comm. Algebra
\textbf{19} (1991), 2081--2095.

\bibitem[C]{Cohnn} P.~M.~Cohn, ``Skew Field Constructions'', London
Math. Soc. Lecture
Notes Series {\bf 27},
Cambridge University Press, Cambridge, 1977.

\bibitem[DF]{DF} L. Diracca and A.
Facchini, \emph{Uniqueness of
monogeny classes for uniform
objects in abelian categories},
J. Pure Appl.
Algebra~\textbf{172} (2002),
183--191.

\bibitem[F1]{libro} A. Facchini,
``Module Theory. Endomorphism
rings and direct sum
decompositions in some classes of
modules", Progress in Math. {\bf
167}, Birkh\"auser Verlag, Basel,
1998.

\bibitem[F2]{JA2002}  A. Facchini,
\emph{Direct sum decompositions
of modules, semilocal
endomorphism rings, and Krull
monoids,} J. Algebra \textbf{256}
(1) (2002), 280--307.

\bibitem[FH]{dolors2} A. Facchini and
D. Herbera, {\em Projective modules over semilocal rings\/}, in
``Algebra and its Applications", D. V. Huynh, S. K. Jain, S. R.
L\'opez-Permouth eds., Contemporary Math. {\bf 259}, Amer. Math.
Soc., Providence, 2000, pp. 181-198.


\bibitem[GO]{gabrieloberst} P.
Gabriel and U. Oberst,
\emph{Spektralkategorien und regul\"are Ringe im Von-Neumann\-schen
Sinn}, Math. Zeitschr.~\textbf{82} (1966), 389--395.

\bibitem[HS]{herberasham} D. Herbera
and A. Shamsuddin,
\emph{Modules with semi-local endomorphism ring}, Proc. Amer.
Math. Soc.~\textbf{123} (1995), 3593--3600.

\bibitem[JL]{jensenlenzing}
C.~U.~Jensen and H.~Lenzing,
``Model Theoretic Algebra",
Algebra, Logic and Applications
Series
\textbf{2}, Gordon and Breach Science Publishers, New York, 1989.

\bibitem[M]{matlis} E. Matlis,
\emph{Some properties of noetherian domains of dimension one},
Canad. J. Math. \textbf{13} (1961), 569--586.

\bibitem[P]{Popescu} N. Popescu,
''Abelian Categories with Applications to Rings and Modules",
Academic Press, London \& New York, 1973.

\bibitem[R]{rowen} L.~H.~Rowen,
\emph{Finitely presented modules over semiperfect rings,} Proc.
Amer. Math. Soc. \textbf{97} (1986), 1--7.

\bibitem[Sc]{schofield}
A.~H.~Schofield,
``Representations of rings over
skew fields", London Math. Soc. Lecture
Notes Series
\textbf{92}, Cambridge University Press, Cambridge, 1985.

\bibitem[St]{stenstrom} B.
Stenstr\"om, ``Rings of
Quotiens'', Springer-Verlag, New
York, 1975.

\bibitem[V1]{Vasconcelos} W. V.
Vasconcelos, \emph{ On finitely
generated flat modules}, Trans.
Amer. Math. Soc.~\textbf{138}
(1969), 505--512.

\bibitem[V2]{Vasc} W. V.
Vasconcelos, \emph{Injective endomorphisms of finitely
generated flat modules}, Proc.
Amer. Math. Soc.~\textbf{25}
(1970), 900--901.

\bibitem[Wa]{warfield} R. B. Warfield
Jr., \emph{Cancellation of
modules and groups and stable
range of endomorphism rings},
Pacific J. Math.~\textbf{91} (2)
(1980), 457--485.

\bibitem[Wi]{Wiegand} R. Wiegand,
{\em Direct-sum decompositions
over local rings,} J.
Algebra~\textbf{240} (2001),
83--97.~\end{thebibliography}
\end{document}